\documentclass{article}
\usepackage{amsmath} 
\usepackage{amssymb}
\usepackage[english]{babel}
\usepackage[latin1]{inputenc}
\usepackage{times}
\usepackage[T1]{fontenc}
\usepackage{graphicx}
\newtheorem{corollary}{Corollary}
		\newtheorem{theorem}{Theorem}

\newcommand  \F {{\mathbb F}}

\usepackage{color}

\begin{document}

\begin{center}
{\Large A note on Hall's sextic residue sequence:\\ correlation measure of order $k$ and related measures of pseudorandomness}\\
\end{center}

\begin{center}
  Hassan Aly$^1$ and Arne Winterhof$^2$
\end{center}
\begin{center}
$^1$ Department of Mathematics, Faculty of Science, Cairo University,\\
P.O. Box 12613, Giza, Egypt\\
E-mail: hassan@sci.cu.edu.eg\\
and\\
Department of Computer Science and Information, College of Science at
Al-Zulfi, Majmaah University, Saudi-Arabia\\

$^2$  Johann Radon Institute for Computational and Applied Mathematics,\\
Austrian Academy of Sciences, Altenberger Str. 69, 4040 Linz, Austria\\
E-mail: arne.winterhof@oeaw.ac.at
\end{center}

\begin{abstract}
 It is known that Hall's sextic residue sequence has some desirable features of pseudorandomness: an ideal two-level autocorrelation and linear complexity
 of the order of magnitude of its period $p$.  
 Here we study its correlation measure of order~$k$ and show that it is, up to a constant depending on $k$ and some logarithmic factor, of order of magnitude $p^{1/2}$, which is close to the expected value for a random sequence of length $p$. Moreover, we derive from this bound a lower bound on the $N$th maximum order complexity of order of magnitude $\log p$, which is the expected order of magnitude for a random sequence of length $p$. 
\end{abstract}



\section{Introduction}

\subsection{Hall's sextic residue sequence}
For a prime $p$ of the form $p=6f+1$,
{\em Hall's sextic residue sequence} ${\cal H}=(h_n)$ of period $p$ is defined as follows:
Let $g$ be a primitive root modulo $p$ and 
\begin{equation}\label{cycl} C_{\ell}=\{g^{6i+\ell} | i=0,1,\cdots, f-1 \},\quad \ell=0,1,\cdots,5,
\end{equation}
be the {\em cyclotomic cosets modulo $p$ of order $6$}.
Then we put 
\begin{equation}
\label{1}
h_n= \left\{\begin{array}{ll} 1&  \mbox{if }  n\bmod p \in C_0 \cup C_1 \cup C_3,\\
                   0 &  \mbox{otherwise},         
    \end{array}\right. \quad n=0,1,\ldots
\end{equation}

Hall's sextic sequence has several desirable features of pseudorandomness:
\begin{itemize}
\item It has a small out-of-phase autocorrelation. In particular, if $p=4u^2+27$ and $g$ is chosen such that $3\in C_1$, then it has ideal $2$-level
autocorrelation (or equivalently $C_0\cup C_1\cup C_3$ is a difference set), see \cite{H56}.
\item It has large linear complexity of order of magnitude $p$ over $\F_2$, see \cite{KS01}.
\item It has large linear complexity over some other fields, see \cite{ES17,HHL16}.
\item It has large $k$-error linear complexity over $\F_p$ for $k<(p-1)/2$, see \cite{AMW07}.
\item It has maximum $2$-adic complexity, see \cite{XQL14} and also \cite{H14} for a short proof.
\end{itemize}

All these features of pseudorandomness consider a full period of the sequence. However, for cryptographic applications usually only a part of a period
of the sequence is used. In this paper we deal with some aperiodic measures of pseudorandomness, the correlation measure of order $k$ and the $N$th maximum order complexity. 

\subsection{Correlation measure of order $k$}

The {\it correlation measure of order $k$} of a binary sequence $(s_n)$ of length $N$ is defined as 
$$ C_k(s_n)= \max_{M,D} \left | \sum_{n=0}^{M-1} {(-1)^{s_{n+d_1}+ s_{n+d_2}+ \cdots+ s_{n+d_k}}} \right |,$$
where the maximum is taken over all $D=(d_1,d_2,\cdots, d_k)$  with non-negative integers $d_1 < d_2 < \cdots < d_k$ and $M$ 
such that $M-1+d_k \le N-1$. 
This measure of pseudorandomness was introduced by Mauduit and S\'ark\"ozy in \cite{MS97}. 

A sequence ${\cal S}=(s_n)$ is considered a good pseudorandom sequence if the value of $C_k({\cal S})$ (at least for small $k$) is small in terms of $N$.

The main result of this paper is the following upper bound on the correlation measure of order $k$ of Hall's sextic sequence.
\begin{theorem}
\label{cor}

Let ${\cal H}=\{h_0,h_1,\cdots,h_{p-1}\}$ be a period of Hall's sextic sequence defined by $(\ref{1})$. Then the correlation measure of order $k$ of ${\cal H}$ satisfies $$C_k({\cal H}) =O\left(\left ( \frac{14}{3}\right )^k k p^{1/2} \log p \right).$$
\end{theorem} 
Here we used the notation $A=O(B)$ if $A\le c B$ for a positive absolute constant $c$. 	

We will prove Theorem~\ref{cor} in Section~\ref{proof}.
    
By \cite{AKMMR07} for a sequence ${\cal S}$ of length $N$ with very high probability the correlation measure $C_k({\cal S})$ is up to a constant depending on $k$ of order of magnitude $\sqrt{N\log N}$. 
In this sense Hall's sextic sequence behaves (almost) like a random sequence.

Moreover, Theorem~\ref{cor} implies bounds on two other measures of pseudorandomness for parts of a period of a sequence, the $N$th maximum order complexity
and the $N$th linear complexity.

Note that bounds on the correlation measure of order $k$ of characteristic sequences of {\em consecutive} unions of cyclotomic classes of order $m$ were studied in \cite{G06}. 
However, Hall's sextic residue sequence is not covered by the construction of \cite{G06} 
and our approach is slightly different.

More precisely, $\cite{G06}$ also deals with the case that the running index $n$ is substituted by~$f(n)$ for a polynomial~$f(X)$ over $\F_p$ of degree $k$ which 
satisfies certain conditions.
However, the result is trivial if $k>p^{1/2}$ but the polynomial used to define Hall's sequence has larger degree $k\ge (p+5)/6$:

Let $C_i'$ be the $i$th cyclotomic class of order $6$ with respect to the primitive element $g^{-1}$ and note that $C_{6-i}=C_i'$ for $i=0,1,\ldots,5$.
Then let $f(n)$ be the mapping which interchanges $C_2$ and $C_3$, that is, 
$$f(n)=\left\{\begin{array}{cc} gn, & n\in C_2,\\
        g^{-1}n, & n\in C_3,\\
        n, & \mbox{otherwise}.
       \end{array}\right.$$
Then we have
$$h_n=\left\{\begin{array}{cc} 0, & 1\le (\mbox{ind}_{g^{-1}}(f(n))\bmod 6)\le 3,\\ 1, & \mbox{otherwise},\end{array}\right. n=1,\ldots,p-1,$$
which is of the form of sequences which are studied in \cite{G06}. (Note that actually \cite{G06} deals with the sequence $((-1)^{h_n})_{n=0}^{p-1}$ over $\{-1,1\}$.) 
Obviously, $f$ cannot be represented by a polynomial $f(X)$ of degree one and is of the form
$$f(X)=X\sum_{i=0}^5 A_i X^{i(p-1)/6}$$
by \cite[Theorem~1]{niwi}.
Hence, $\deg f \ge (p+5)/6$ and $\cite{G06}$ does not give a nontrivial bound.

\subsection{$N$th Maximum order complexity}

The $N$th {\it maximum order complexity} $M({\cal S},N)$ of a binary sequence ${\cal S}=(s_n)$ is the smallest positive integer $M$
such that there is a polynomial $f(x_1,\ldots,x_M)\in \F_2[x_1,\ldots,x_M]$
with
$$s_{i+M}=f(s_i,s_{i+1},\ldots,s_{i+M-1}),\quad 0\le i\le N-M-1,$$
see \cite{ja89,ja91,nixi14}. 

By \cite{IW17} for any binary sequence ${\cal S}$ we have the following relation between the maximum order complexity and the correlation measure of order $k$:
$$M({\cal S},N) \ge N- 2^{M({\cal S},N)+1} \max_{1\le k \le M({\cal S},N)+1} C_k({\cal S},N),\quad N\ge 1.$$
Combining this relation and Theorem~\ref{cor} we get the following lower bound on the maximum order complexity of Hall's sextic residue sequence.
\begin{corollary}
\label{MH}
For the $N$th maximum order complexity of Hall's sextic residue sequence we have
$$M({\cal H},N) = \Omega \left( \log \left(\frac{\min\{N,p\}}{p^{1/2}\log^2p}\right)\right).$$ 
\end{corollary}
Here $A=\Omega(B)$ is equivalent to $B=O(A)$.

The expected value of the $N$th maximum-order complexity is of order of magnitude $\log N$, see~\cite{ja89} as well as \cite{nixi14} (Remark 4) and references therein. Hence, the $N$th maximum order complexity of Hall's sextic residue sequence for any $N$ with $p^{1/2}\log^3 p\le N\le p$ is at least of this desired order of magnitude.

\subsection{$N$th linear complexity}

For $N\ge 1$ the {\it $N$th linear complexity} $L({\cal S},N)$ over $\F_2$ of a binary sequence ${\cal S}=(s_n)$ is the shortest Length $L$ of a linear recurrence relation over $\F_2$
$$ s_{n+L}=c_{L-1}s_{n+L-1}+ \cdots + c_0s_n, 0 \le n \le N-L-1,$$
which is satisfied by the first $N$ sequence elements. (We use the convention that $L({\cal S},N)=0$ if the first $N$ terms of the sequence are all 0 and $L({\cal S},N)= N$ if $ S_0=s_1=\cdots =s_{N-2}=0$ and $s_{N-1}=1$.) 

\cite[Theorem~1]{BW06} gives a lower bound on the $N$th linear complexity of ${\cal S}$ in terms of the correlation measure of order $k$: 
$$ L({\cal S}, N) \geq N- \max_{1 \leq k \le L({\cal S}, N)+1} C_k({\cal S}),\quad N\ge 1.$$ 
Combining this bound with Theorem~\ref{cor} we get for Hall's sextic residue sequence ${\cal H}$ 
$$ L({\cal H},N) = \Omega \left( \log \left(\frac{\min\{N,p\}}{p^{1/2}\log^2p}\right)\right).$$ 
The implied constant can be explicitly calculated.
This result can be even obtained in a simpler way from Corollary~\ref{MH} observing that
$$L({\cal H}, N)\ge M({\cal H},N)$$
with a slightly weaker implied constant. 

Some experiments indicate that $L({\cal H},N)$ is much larger and we consider it an interesting open problem to improve this lower bound
on $L({\cal H},N)$.

Note that there is another figure of merit closely related to the correlation measure of order $k$, the {\em arithmetic autocorrelation}. Rather moderate bounds on the arithmetic autocorrelation of Hall's sextic residue sequence can be obtained combining Theorem~\ref{cor} and \cite{HMW17}.




\section{Proof of Theorem~\ref{cor}}
\label{proof}

Let  $\eta$ be a multiplicative character of $\F_p$ of order $3$ and put
$$\delta_1(n)=\frac{1+\eta(n)+\eta^2(n)}{3}, \quad n\in \F_p^*.$$ 
Then we have  
$$
\delta_1(n)=\left\{\begin{array}{ll} 1 &  \mbox{  if }  n \in C_0 \cup C_3,\\
                   0 & \mbox{  if }  n \in C_1 \cup C_2 \cup C_4 \cup C_5.                  
    \end{array}\right.
$$
Let $\chi$ be a multiplicative character of $\F_p$ of order $6$
and put 
$$\omega = \chi(g),$$
where $g$ is the primitive root modulo $p$ fixed for the definition of the cyclotomic cosets in $(\ref{cycl})$
and 
$$\delta_2(n)  = \frac{1+\omega^{-1}\chi(n)+ \omega^{-2}\chi^2(n) +\cdots \omega^{-5}\chi^5(n)}{6}.$$
Then we have
\begin{equation}\nonumber
\begin{split}
\delta_2(n) 
 & = \left\{\begin{array}{ll} 1,&  \mbox{  if }  n \in C_1,\\
                   0, & \mbox{  if } n \in C_0 \cup C_2 \cup C_3 \cup C_4 \cup C_5,               
    \end{array}\right. \quad n\in \F_p^*.
\end{split}
\end{equation}
Now $\delta(n)$ defined by
$$\delta(n)=\delta_1(n)+\delta_2(n),\quad n\in \F_p^*,$$
is the characteristic function of $C_0\cup C_1\cup C_3\subset \F_p^*$.
Hence, Hall's sextic sequence ${\cal H}=(h_n)$ defined by $(\ref{1})$ satisfies
$$
h_n= \delta(n),\quad n=1,2,\ldots,p-1, 
$$		
and we have
\begin{equation*}
\begin{split}
 (-1)^{h_n} & = 1-2\delta(n) \\
 & = \frac{-2}{3}\left(\eta(n)+\eta^2(n)\right)-\frac{1}{3}\left(\omega^{-1}\chi(n)+\omega^{-2}\chi^2(n)+\cdots+ \omega^{-5}\chi^5(n)\right)
\end{split}
\end{equation*}
for $n=1,2,\ldots,p-1$.
Note that $\eta\in \{\chi^2,\chi^4\}$.
Now
$$\left|\sum_{n=1}^{M-1} (-1)^{h_{n+d_1}+\ldots+h_{n+d_k}}\right|$$
can be estimated by $7^k$ sums of the form
$$\left(\frac{2}{3}\right)^k\left|\sum_{n=1}^{M-1} \chi((n+d_1)^{m_1}\cdots (n+d_k)^{m_k})\right|$$
with $1\le m_1,\ldots,m_k\le 5$.
A variant of Weil's theorem for incomplete character sums, see for example \cite[Lemma~3.4]{S03}, gives the bound
$$O(kp^{1/2}\log p)$$
for the absolute value of the inner character sums.
Collecting everything, Theorem~\ref{cor} follows.~\hfill $\Box$
%

\section{Legendre sequence and Ding-Helleseth-Lam sequence}

There are several other cyclotomic sequences with similar features of pseudorandomness. We mention only the Legendre sequence and the Ding-Helleseth-Lam sequence.

\subsection{Legendre sequence}
Similar results are known for the Legendre sequence ${\cal L}=(\ell_n)$ of prime period $p>2$ which is
the characteristic sequence of the quadratic residues modulo $p$, that is, a cyclotomic sequence of order $2$:
\begin{itemize}
\item It has small out-of-phase autocorrelation and in particular ideal $2$-level autocorrelation if $p\equiv 3\bmod 4$, see \cite{P33}. 
\item It has large linear complexity of order of magnitude $p$, see \cite{DHS98,T64}. 
\item It has $k$th error linear complexity over $\F_p$ of order of magnitude $p$ for $k<(p-1)/2$, see \cite{AW06}. 
\item Its $2$-adic expansion is maximal, see \cite{XQL14,H14,HW18}. 
\item It has small correlation measure of order $k$ of order of magnitude $kp^{1/2}\log p$, see \cite{MS97}.
\item Its $N$th maximum order complexity is at least $\log(\min\{N,p\}/p^{1/2})+O(\log\log p)$, see \cite{IW17}.
\item Its $N$th linear complexity is at least of order of magnitude $\min\{N,p\}/(p^{1/2}\log p)$, see \cite[Theorem 9.2]{S03}
or combine \cite{MS97} and \cite{BW06}.
\item Its arithmetic autocorrelation is moderately small, see \cite{HW17}.
\end{itemize}

\subsection{Ding-Helleseth-Lam sequence}

Ding et al.\ \cite{DHL99} introduced a cyclotomic generator of order $4$. Let $C_0=\{x^4: x\in \F_p^*\}$ be the subgroup of $\F_p^*$ of 
bi-squares and $C_1=gC_0$ for a primitive root $g$ modulo $p$. Then the {\em Ding-Helleseth-Lam sequence} ${\cal D}=(d_n)$ is the characteristic sequence of $C_0\cup C_1$.
\begin{itemize}
\item Its out-of-phase autocorrelation is small and it has optimum three-level autocorrelation (or equivalently $C_0\cup C_1$ is an almost difference set)
if $p=x^2+4$ with $x\equiv 1\bmod $, see \cite[Theorem~4]{DHL99}.
\item It has linear complexity close to its period, see \cite{DHL00}, if $2\in D_0$.
\item It has maximum $2$-adic complexity, see \cite{XQL14,H14}.
\item For $k<(p-1)/2$ it has $k$-error linear complexity over $\F_p$ of order of magnitude $p$, see \cite{AMW07}. 
\item Its correlation measure of order $k$ is estimated in \cite{G06}.
\item Its $N$th maximum order complexity can be lower bounded by combining \cite{IW17} and~\cite{G06}.
\end{itemize}



It is possible to extend the results of this paper to cyclotomic sequences of higher order $m$. 
For the special case of characteristic sequences of the union of consecutive cyclotomic classes see \cite{G06}.
However, for odd $m$ the correlation measure of order $k$ can be very large, see \cite{G06}, and for even $m$ the implied constant in the bound on the correlation measure of order $k$ either depends on $m$ 
(non-consecutive case) or we have an additional factor of $(\log p)^k$ (consecutive case \cite{G06}). Hence, the Hall sequence, the Legendre sequence and the Ding-Helleseth-Lam sequence are certainly the
most attractive candidates for cryptographic applications.

Moreover, there are several other sequence constructions using characters of finite fields or elliptic curves over finite fields, see for example \cite{G13} and references therein.

\section*{Acknowledgment}
This paper was partly written during a pleasant visit of the first author at RICAM. He wishes to thank for hospitality and financial support.
The second author is supported by the Austrian Science Fund FWF Project P~30405-N32.

\end{document}